\input amssymb.sty
\documentclass{amsproc}
\usepackage{graphicx}
\newtheorem{theorem}{Th\'eor\`eme}[section]

\theoremstyle{definition}

\newtheorem{proposition}[theorem]{Proposition}
\theoremstyle{remark}

\newcommand{\og}{\leavevmode\raise.3ex\hbox{$\scriptscriptstyle\langle\!\langle$}}
\newcommand{\fg}{\leavevmode\raise.3ex\hbox{$\scriptscriptstyle\rangle\!\rangle$}}
\numberwithin{equation}{section}

\begin{document}

\title{Matrix valued Brownian motion and a paper by  P\'olya}

\author{Philippe Biane}
\address{CNRS, Laboratoire d'Informatique 
 Institut Gaspard Monge 
 Universit\'e Paris-Est 
 5 bd Descartes, Champs-sur-Marne 
 F-77454 Marne-la-Vall\'ee cedex 2
}
\email{Philippe.Biane@univ-mlv.fr}

\maketitle

\section{Introduction}
\label{sec:1}
This paper has two parts which are largely independent.
In the first one I recall some known facts on matrix valued Brownian motion,
which are not so easily found in this form in the literature.
I will study three types of matrices, namely Hermitian matrices,
complex invertible matrices, and unitary matrices, and try to give a precise
description of the motion of eigenvalues (or singular values) in each case.
In the second part, I give a new look at an old paper of 
  G. P\'olya \cite{P}, where he introduces a function close to 
     Riemann's $\xi$ function, and shows that it satisfies 
 Riemann's hypothesis. As put by Marc Kac
in his comments on P\'olya's paper \cite{K}, 
``Although this beautiful paper
takes you within a hair's breadth of Riemann's hypothesis it does not seem to
have inspired much further work and reference to it in the mathematical
litterature are rather scant''.   
My aim is to point out that the function considered by P\'olya
is related in a more subtle way to Riemann's $\xi$ function than it looks at first
sight. Furthermore the nature of this relation is probabilistic, since these
functions have a natural interpretation involving Mellin transforms of
first passage
times for diffusions. By studying infinite divisibility properties of the
distributions of these first passage times, we will see that they  are
generalized gamma convolutions, whose mixing measures are related to the
considerations in the first part of this note.

\section{Matrix Brownian motions}
We will study three types of matrix spaces, and in each of these spaces consider
a natural Brownian motion, and show that the motion of eigenvalues (or singular
values) of this Brownian motion has a simple geometric description, using Doob's
transform. 
The following results admit analogues in more general complex symmetric spaces,
but for the sake of simplicity, discussion will be restricted to 
type  $A$ symmetric spaces. Actually the interesting case for us in the second
part will be the simplest one, of rank one, but I think that
this almost trivial case is better understood by putting it in the more general context.
Some references for results in this section are
 \cite{Ba}, \cite{B}, \cite{BPY},
  \cite{D}, \cite{GLT}, \cite{H}, \cite{JO}, \cite{M},  \cite{T1}, \cite{T2}.

\subsection{Hermitian matrices}Consider the space
of $n\times n$ Hermitian matrices, with zero trace,  endowed with the 
 quadratic form $$\langle A,B\rangle=Tr(AB).$$  Let $M(t)$ be a Brownian motion
with values in this space, which is simply a
Gaussian process with covariance $$E[Tr(AM(t))Tr(BM(s))]=Tr(AB)s\wedge t$$
for $A,B$   traceless Hermitian matrices.

 Let  $\lambda_1(t)\geq \lambda_2(t)\geq\ldots\geq
\lambda_n(t)$ be the eigenvalues of $M(t)$; they 
perform a stochastic process with
values in 
the Weyl chamber $${\mathcal C}=\{(x_1,\ldots, x_n)\in {\bf R}^n\,|
\,x_1\geq x_2\geq \ldots
\geq x_n\}\cap H_n$$
where
$$H_n=\Bigl\{\,(x_1,\ldots, x_n)\in {\bf R}^n\,\Bigm|\,\sum_{i=1}^nx_i=0\,\Bigr\}\;.$$
Let $p_t^0$ be the transition probability 
semi-group  of Brownian motion   killed at the boundary of the
cone ${\mathcal C}$. This cone is a   fundamental domain for the action 
of the symmetric
group $S_n$, which acts by  permutation of coordinates on 
${\bf R}^n$. Using the reflexion principle, one shows easily that 
  $$p^0_t(x,y)=\sum_{\sigma\in
S_n}\epsilon(\sigma)p_t(x,\sigma(y))\qquad x,y\in {\mathcal C}$$
where $p_t(x,y)=(2\pi t)^{-(n-1)/2}e^{-|x-y|^2/2t}$ and
 $\epsilon(\sigma)$ is the  signature of $\sigma$.
Let $h$ be the function  
$$h(x)=\prod_{i>j}(x_i-x_j).$$
\begin{proposition}
The function   $h$ is the unique  (up to a positive multiplicative
  constant)
  positive harmonic function for the  semigroup $p_t^0$, 
on the   cone ${\mathcal C}$, which vanishes on the boundary.
\end{proposition}
 The harmonic function $h$ corresponds to the unique   point at infinity
 in the Martin  compactification of
$\mathcal
C$. Consider now the Doob's transform 
of $p_t^0$, which is the semigroup given by
$$q_t(x,y)=\frac{h(y)}{h(x)}p^0_t(x,y).$$ 
It is a diffusion semigroup on $\mathcal C$ with infinitesimal generator
$$\frac{1}{2}\Delta+\langle \nabla \log h,\nabla \,\cdot\,\rangle\;.$$
\begin{proposition}
The eigenvalue process of a traceless Hermitian Brownian
motion is a Markov  diffusion process in the cone 
$\mathcal C$, with  semigroup $q_t$. 
\end{proposition}
We can summarize the last proposition by saying that
the eigenvalue process is a Brownian motion in 
$\mathcal C$, conditioned 
(in Doob's sense) to exit the cone at infinity.
\subsection{The group $SL_n({\bf C})$}
This is the  group of complex invertible   matrices 
of size  $n\times n$, with
determinant 1. Its Lie algebra is the space 
${\mathfrak s\mathfrak l}_n(\bf C)$
of complex traceless  matrices. 
Consider the Hermitian form
 $$\langle A,B\rangle=Tr(AB^*)$$
 on ${\mathfrak s\mathfrak l}_n(\bf C)$
which is invariant by left and right action of the 
 unitary subgroup $SU(n)$. This Hermitian form  determines a
unique  Brownian motion  with  values in
  ${\mathfrak s\mathfrak l}_n(\bf C)$.
 The Brownian motion $g_t$,  on  
 $SL_n(\bf C)$, is the   stochastic exponential of this 
Brownian motion, solution to the Stratonovich stochastic 
differential equation 
$$dg_t=g_tdw_t$$
where  $w_t$ is a
 Brownian motion in
${\mathfrak s\mathfrak l}_n(\bf C)$.

There are two remarkable decompositions of 
$SL_n(\bf C)$,  the Iwasawa and Cartan decompositions.
The first one is $SL_n({\bf C})=NAK$ where 
 $K$ is the compact group
 $SU(n)$,  $A$ is the group of diagonal matrices 
with positive coefficients, and determinant one, and $N$
is the   nilpotent  group of upper triangular  matrices 
with 1's on the diagonal. 
Each matrix of  $SL_n(\bf C)$ has a unique 
decomposition as a product $g=nak$ of elements of the three subgroups
 $N,A,K$.
This can be easily inferred from the Gram-Schmidt orthogonalization process.
If $g_t$ is a  Brownian motion in  $SL_n({\bf C})$, one can consider its
components 
 $n_t, a_t, k_t$. 
In  particular, denoting by $(e^{w_1(t)},\ldots, e^{w_n(t)})$ the diagonal
 components of
$a_t$ the following holds (cf \cite{T1}).
\begin{proposition}
The process $\bigl(w_1(t),\ldots, w_n(t)\bigr)$ is a  Brownian motion with a 
drift
$\rho=(-n+1,-n+3,\ldots, n-1)$ in the subspace $H_n$.
\end{proposition}

The other decomposition is the Cartan decomposition 
 $SL_n({\bf C})=KA^+K$, where 
$A^+$ is the part of $A$ consisting of
 matrices with positive nonincreasing  coefficients along the diagonal.
In order to get the Cartan decomposition of a matrix 
$g\in SL_n(\bf C)$, take its  polar decomposition $g=ru$ with
$r$  positive Hermitian,
and $u$ unitary, then diagonalize $r$ which yields 
$g=vav'$ with $v$ and $v'$ unitary and $a$ diagonal, with positive real
 coefficients
which can be put in nonincreasing order along the diagonal.
These  coefficients
are the singular values of 
 the matrix $g$. This  decomposition is not
unique since the diagonal subgroup of 
$SU(n)$ commutes with   $A$, but the singular values are uniquely defined.
Call  $(e^{a_1(t)},\ldots ,e^{a_n(t)})$ the singular values of 
the Brownian  motion $g_t$, with  
$a_1\geq a_2\geq \ldots\geq a_n$. They form a process with values
 in the  cone $\mathcal C$. Let us mention that this stochastic process can also
 be interpreted as the radial part of a Brownian motion with values in the 
 symmetric space  $SL_n({\bf C})/SU(n)$.
 We will now give for the motion of singular values  a similar 
 description as the one of eigenvalues of the Hermitian Brownian motion.
 For this, consider a  Brownian motion in  $H_n$, 
with drift $\rho$, killed at the exit of 
the cone $\mathcal C$. This process has a semigroup
 given by
$$p^{0,\rho}_t(x,y)=e^{\langle \rho,y-x\rangle-
t\langle \rho,\rho\rangle/2}
p^0_t(x,y)\;.$$
\begin{proposition}
The function 
$$h^\rho(y)=\prod_{i>j} (1-e^{2(y_j-y_i)})$$ is a positive harmonic function
for the semigroup $p^{0,\rho}_t$,
in the cone $\mathcal C$, and vanishes at the boundary of the cone.
\end{proposition}  
It is not true that this function is the unique
positive harmonic function on the cone; indeed the
Martin boundary at infinity is now much larger  and contains 
a point for each  direction inside the cone, see
  \cite{GLT}.
The Doob-transformed semigroup
$$q^\rho_t(x,y)=\frac{h^\rho(y)}{h^\rho(x)}p^{0,\rho}_t(x,y)$$
is a Markov    diffusion semigroup
in the  cone $\mathcal C$, with infinitesimal generator 
$$\frac{1}{2}\Delta+\langle \rho,\cdot\rangle+\langle \nabla\log h^\rho,
\nabla\cdot\rangle.$$
Note that   it can also be expressed as
$$\frac{1}{2}\Delta+\langle \nabla\log \tilde h^\rho,
\nabla\cdot\rangle$$
with  
$$\tilde h^\rho(y)=\prod_{i>j} \sinh(y_j-y_i)$$
(see\cite{JO}).
\begin{proposition}
The logarithms of the singular values of a Brownian motion  in 
$SL_n(\bf C)$ perform a diffusion process
 in the cone $\mathcal C$,
 with semigroup
$q^\rho_t$. 
\end{proposition}  
As in the preceding case, we can summarize by saying that
the process of singular
values is a Brownian  motion with drift $\rho$
 in the cone 
$\mathcal C$, conditioned 
(in Doob's sense) to exit the cone at infinity, in the direction $\rho$.
\subsection{Unitary matrices }
The Brownian motion with values in $SU(n)$ is obtained by taking the stochastic
exponential of 
a Brownian motion in the Lie algebra of traceless 
anti-Hermitian matrices, endowed with the Hermitian form
 $$\langle A,B\rangle=-Tr(AB).$$
Let 
$e^{i\theta_1},\ldots, e^{i\theta_n}$ be the eigenvalues of a matrix in $SU(n)$,
which can be chosen so that 
$\sum_i\theta_i=0$, and
  $\theta_1\geq \theta_2\geq\ldots\geq\theta_n$, $\theta_1-\theta_n\leq 2\pi$.
   These  conditions
determine a simplex $\Delta_n$ in $H_n$, which is a  fundamental domain
 for the action of the affine
 Weyl group on  $H_n$. Recall that the affine 
   Weyl group  $\tilde W$ is the semidirect product of the 
   symmetric group
$S_n$, which acts by permutation of coordinates in  $H_n$,
 and of the group of 
translations by elements of the lattice
 $(2\pi{\bf Z})^n\cap H_n$.

One can use the reflexion principle again to compute the semigroup
of Brownian moton in this simplex killed at the boundary.
One gets an alternating sum over the elements of $\tilde W$, 
$$p^0_t(\theta,\xi)=\sum_{w\in\tilde W}\epsilon(w)p_t(\theta,w(\xi)).$$
The infinitesimal generator
is 1/2 $\times$ the Laplacian in the simplex, with Dirichlet boundary
conditions. It is well known that this operator has a compact resolvent, and its
eigenvalue with smallest module is simple, with an eigenfunction which can be
chosen positive.
Consider the function
$$h^u(\theta)=\prod_{j> k}(e^{i\theta_j}-e^{i\theta_k}).$$
\begin{proposition}
The function $h^u$ is positive  inside the simplex $\Delta_n$, it vanishes on
the boundary, and it is the eigenfunction
corresponding to the  Dirichlet  eigenvalue 
 with smallest module 
on $\Delta_n$. This eigenvalue is
$\lambda=(n-n^3)/6$.
\end{proposition}
The Doob-transformed semigroup
$$q^u_t(x,y)=\frac{h^u(y)}{h^u(x)}e^{-\lambda t}p^0_t(x,y)$$
is a Markov  diffusion semigroup in $\Delta_n$, with infinitesimal generator
$$\frac{1}{2}\Delta+\langle \nabla\log h^u,\nabla\cdot\rangle-\lambda.$$
\begin{proposition}
The process of eigenvalues of a unitary Brownian motion 
is a diffusion with values in  
$\Delta_n$ with probability transition semigroup $q^u_t$.
\end{proposition}
Again  a good summary of this situation is that the motion of eigenvalues is
that of a Brownian motion in the simplex $\Delta_n$ conditioned to stay forever 
in this
simplex.
\subsection{The case of rank 1}\label{rg1}
In the next section we will need the simplest case, that of 
$2\times 2$ matrices.
Consider first the case of Hermitian matrices. The process of
eigenvalues is essentially a Bessel process of
 dimension 3, with infinitesimal generator
$$\frac{1}{2}\frac{d^2}{dx^2}+\frac{1}{x}\frac{d}{dx}\qquad \hbox{on}\quad
]0,+\infty[,$$
obtained from Brownian motion killed at zero, of infinitesimal generator 
$$\frac{1}{2}\frac{d^2}{dx^2}$$with Dirichlet boundary condition at 0, 
by a Doob transform with the positive
harmonic function $h(x)=x$

In the case of the  group $SL_2({\bf C})$, or the symmetric space 
$SL_2({\bf C})/SU(2)$, which is the hyperbolic 
space of dimension 3, the radial
process has infinitesimal generator
$$\frac{1}{2}\frac{d^2}{dx^2}+\coth x\frac{d}{dx}$$
obtained from Brownian motion with a drift
$$\frac{1}{2}\frac{d^2}{dx^2}+\frac{d}{dx}$$
with  Dirichlet boundary  condition at 0, by a Doob's transform 
with the function $1-e^{-2x}$.

Finally the last case is Brownian motion in
 $SU(2)$, where the eigenvalue process takes values in 
 $[0,\pi]$ and has an infinitesimal generator 
$$\frac{1}{2}\frac{d^2}{dx^2}+\cot x\frac{d}{dx}$$
obtained by a Doob transform at the bottom of the spectrum
from
$$\frac{1}{2}\frac{d^2}{dx^2}$$ on $[0,\pi]$ with Dirichlet boundary conditions
at $0$ and $\pi$, by the function
 $\sin(x)$.

For these last two examples, we shall write a spectral decomposition of the 
generator $L_i,i=1,2$, of the form
\begin{equation}\label{spectral}
f(x)=\int \Phi^i_{\lambda}(x)\Bigl[\int\Phi^i_{\lambda}
(y)f(y)\,dm_i(y)\Bigr]\,d\nu_i(\lambda)\qquad i=1,2
\end{equation}
for every $f\in L^2(m_i)$,
where  $m_i$ is measure for which $L_i$ is selfadjoint in
$L^2(m_i)$,  and the functions  $\Phi^i_\lambda$ are solutions to 
$$L_i\Phi_\lambda^i+\lambda\Phi_\lambda^i=0$$
and $\nu_i$ is a spectral measure for $L_i$.

For $L_1=\frac{1}{2}\frac{d^2}{dx^2}$ on $[0,\pi]$
with Dirichlet boundary conditions,
  $m_1(dx)$ is
 Lebesgue measure on $[0,\pi]$, and  $L_1$ is selfadjoint on 
  $L^2([0,\pi)]$.
 Furthermore $$\Phi_\lambda^1(x)=\sin(\sqrt{2\lambda}x)$$
 and
\begin{equation}\label{nu1}
\nu_1(d\lambda)=\frac{1}{\pi}\sum_{n=1}^\infty \delta_{n^2/2}(d\lambda).
\end{equation}
 For $L_2=\frac{1}{2}\frac{d^2}{dx^2}+\frac{d}{dx}$ on $[0,+\infty[$,
 the measure $m_2(dx)=e^{2x}dx$, and
 $$\Phi^2_\lambda(x)=e^{-x}\sin(\sqrt{2\lambda-1}x)\qquad \lambda>1/2.$$ 
 The spectral measure is
 \begin{equation}\label{nu2}
\nu_2(d\lambda)=\frac{1}{\pi\sqrt{2\lambda-1}}d\lambda\qquad \lambda>1/2
\end{equation}
 on $[1,+\infty[$.
Of course  formulas (\ref{spectral}), (\ref{nu1}), (\ref{nu2})
 are immediate consequences of ordinary 
 Fourier analysis.

 Note that the spectral decompositions, and in particular 
the measures $\nu_i$, 
 depend on the normalisation of the fonctions $\Phi_\lambda$.
 We have made a natural choice, but it  does not coincide with the usual
 normalisation of  Weyl-Titchmarsh-Kodaira theory, see \cite{CL}.

\section{ MacDonald's function and Riemann's $\xi$ function}
\subsection{ P\'olya's paper}  
In his paper \cite{P}, P\'olya starts from Riemann's $\xi$ function
$$\xi(s)=s(s-1)\pi^{-s/2}\Gamma(s/2)\zeta(s)$$
where  $\zeta$ is Riemann's zeta function. Then $\xi$ is an entire function
whose zeros are exactly the nontrivial zeros of  
$\zeta$.

Putting $s=1/2+iz$ yields
$$\xi(z)=2\int_0^\infty\Phi(u)\cos(zu)du$$
with
\begin{equation}\label{theta}\Phi(u)=2\pi 
e^{5u/2}\sum_{n=1}^\infty (2\pi e^{2u}n^2-3)n^2e^{-\pi n^2e^{2u}}
\end{equation}
and the function $\Phi$ is even, as follows from the functional equation for 
Jacobi
 $\theta$ function; furthermore 
$$\Phi(u)\sim 4\pi^2e^{9u/2-\pi e^{2u}}\qquad u\to+\infty$$
so that
$$\Phi(u)\sim 4\pi^2(e^{9u/2}+e^{-9u/2})
e^{-\pi (e^{2u}+e^{-2u})}\qquad u\to\pm\infty.$$
This lead P\'olya to define  a ``falsified'' $\xi$ function
$$\xi^*(z)=8\pi^2\int_0^\infty(e^{9u/2}+e^{-9u/2})
e^{-\pi (e^{2u}+e^{-2u})}\cos(zu)du.$$
The main result of \cite{P} is
\begin{theorem} The function $\xi^*$ is entire, its 
 zeros are  real and  simple. 
Let $N(r)$,
(resp. $N^*(r)$) denote the number of zeros of $\xi(z)$ (resp. $\xi^*(z)$)
 with real part in the interval $[0,r]$, then 
 $N(r)-N^*(r)=O(\log r)$.
\end{theorem}
Recall  that the same assertion about the zeros of  
the function $\xi$ (without the statement about simplicity, beware also that $s=1/2+iz$) is Riemann's hypothesis.
Recall also the well known estimate  
$$N(r)= \frac{r}{2\pi}\log (r/2\pi)-\frac{r}{2\pi}+O(1).$$
P\'olya's results rely on the intermediate study of the  function
$$\mathfrak{G}(z,a)=\int_{-\infty}^{\infty}
 e^{-a(e^u+e^{-u})+zu}du$$ from which  $\xi^*$ is obtained by 
$$\xi^*(z)=2\pi^2(\mathfrak{G}(iz/2-9/4,\pi)+\mathfrak{G}(iz/2+9/4,\pi))$$
P\'olya shows that  $\mathfrak{G}(z,a)$ has only purely imaginary zeros, (as a function of $z$)
and the  number of these zeros with imaginary part in $[0,r]$ 
grows as $\frac{r}{\pi}\log\frac{r}{a}-\frac{r}\pi+O(1)$.
The results on $\xi^*$ are then deduced through a nice lemma  which 
played a role in the history of statistical mechanics (the Lee-Yang theorem on Ising model), as revealed by 
M. Kac \cite{K}. We shall now concentrate on  $\mathfrak{G}(z,a)$. In particular, for
  $a=\pi$,  the function
$\tilde \xi(z)=\mathfrak{G}(iz/2,\pi)$  is another approximation of 
$\xi$ which has many interesting structural properties.
\subsection{ MacDonald functions} 
The function denoted $\mathfrak{G}(z,a)$ by P\'olya is actually a Bessel function.
Indeed, MacDonald's function, also called modified Bessel function (see e.g. \cite{A}), given by
$$K_{\mu}(x)=\int_0^\infty t^{\mu-1}e^{-\frac{x}{ 2}(t+t^{-1})}dt\qquad
x>0,\, \mu\in{\bf C}.$$
satisfies $K_{z}(2x)=\mathfrak{G}(z,x)$. The function $\mathfrak{G}(z,a)$  is therefore
essentially a MacDonald function, as noted by P\'olya.
 MacDonald function is an even function of  $\mu$ and satisfies
\begin{equation}\label{K1}
\frac{2\mu}{ x} K_{\mu}(x)=K_{\mu+1}(x)-K_{\mu-1}(x)
\end{equation}
\begin{equation}\label{K2}
-2\frac{d}{dx}K_{\mu}(x)=K_{\mu+1}(x)+K_{\mu-1}(x)
\end{equation}
The first of these equations is used by 
P\'olya in a very clever way to 
prove that the zeros (in $z$) of 
$\mathfrak{G}(z,x)$ are purely imaginary.

\subsection{Spectral interpretation  of the zeros}
From (\ref{K1}), (\ref{K2})
$$\begin{array}{c} (\frac{\mu}{ x}-\frac{d}{dx})K_{\mu}=K_{\mu+1}\\ \\
(-\frac{\mu}{ x}-\frac{d}{ dx})K_{\mu}=K_{\mu-1}
\end{array}$$
from which one gets
$$\begin{array}{c}
K_{\mu}=(-\frac{\mu+1}{ x}-\frac{d}{ dx})
(\frac{\mu}{ x}-\frac{d}{dx})K_{\mu}\\ \\
=(\frac{d^2}{ dx^2}+\frac{1}{x}\frac{d}{dx}-\frac{\mu^2}{ x^2})K_{\mu}\;.
\end{array}$$
This differential equation will give us a spectral interpretation of the zeros of 
$\mathfrak{G}(z,x)$.
Change  variable by $\psi_{\mu}(x)=K_{\mu}(e^x)$ to get
\begin{equation}\label{STL}
(-\frac{d^2}{ dx^2}+e^{2x})\psi_{\mu}=-\mu^2\psi_{\mu}
\end{equation}
Since $K_{\mu}$ vanishes exponentially at infinity, the spectral theory of Sturm-Liouville operators
on the half-line
(see e.g. \cite{CL}, \cite{LS}) implies that
the squares of the zeros of $\mu\mapsto\psi_{\mu}(y)$ are the eigenvalues of 
$\frac{d^2}{ dx^2}-e^{2x}$ on the interval
$[y,+\infty[$ with the Dirichlet boundary condition at $y$, the functions 
$\psi_\mu$ being the eigenfunctions. Since this operator 
is selfadjoint and negative
 the zeros are  purely imaginary, and are simple.

This spectral interpretation of the zeros of MacDonald function is well known
 \cite{T}, I do not know why P\'olya does not mention it.

\subsection{$H=xp$}
Equation (\ref{STL}) can be put into Dirac's form, indeed the equations
$$\begin{array}{rl} \left( \frac{d}{dx}+\frac{1}{ 2}+e^x\right)f&=\gamma g\\
\\
\left(-\frac{d}{dx}+\frac{1}{ 2}+e^x\right)g&=\gamma f\end{array}$$
imply
$$\left(-\frac{d^2}{dx^2}+e^{2x}\right)f=(\gamma^2-\frac{1}{4})f.$$
Using the change of variables $u=e^x$,
  we get
$$\begin{array}{rl} \left( u\frac{ d}{du}+\frac{1}{2}+u\right)f&=\gamma g\\
\\
\left(-u\frac{d}{du}+\frac{1}{2}+u\right)g&=\gamma f.\end{array}$$
Remark that this  Dirac system yields a perturbation of the 
 Hamiltonian $H=xp$
considered by  Berry  et Keating \cite{BK}, in relation with
 Riemann's zeta function.
\subsection{Asymptotics of the zeros}
General results on Sturm-Liouville operators allow one to recover the asymptotic
behaviour of the spectrum, thanks to a semiclassical analysis, see e.g.
\cite{LS}.
One can get a more precise result using the integral representation of 
 $K_{i\mu}$.
P\'olya gives the asymptotic estimate
$$K_{x+iy}(2a)=
\frac{1}{ \sqrt{2\pi y}}
e^{-\frac{\pi}{2}y+
i\frac{\pi}{2}x}
\left[\left(
\frac{y}{ a}\right)^xe^{i\Phi}+
\left(\frac{y}{a}\right)^{-x}
e^{-i\Phi}\right]+
O(e^{-\frac{\pi}{ 2}y}y^{|x|-3/2})$$
in the strip $|x|\leq 1$ uniformly as $y\to \infty$, where
$$\Phi=y\log\frac{y}{ a}-y-\frac{\pi}{ 4}\;.$$
This estimate can be obtained by the stationary phase method, writing
$$K_{z}(2a)=\int_{-\infty}^\infty e^{zt-2a\cosh(t)}dt.$$ 
Making a contour deformation we get
$$\begin{array}{rl}
K_{z}(2a)=&\int_{-\infty}^{-A} e^{zt
-2a\cosh(t)}dt+i\int_0^{\pi/2}e^{z(-A+it)-2a\cosh(-A+it)}dt\\&\quad
+\int_{-A}^A
e^{z(t+i\frac{\pi}{ 2})-2ai\sinh(t)}dt -
i\int_0^{\pi/2}e^{z(A-it+i\frac{\pi}{ 2})-2a\cosh(A-it+i\frac{\pi}{ 2})}dt
\\&\qquad+\int_A^\infty e^{zt-2a\cosh(t)}dt
\end{array}$$ 
and P\'olya's estimate can be obtained by standard methods, which give also
estimates for the derivatives of
  MacDonald's function.
Finally the zeros of $y\to K_{iy}(2a)$ behave like the 
 solutions to 
$$y\log\frac{y}{
a}-y-\frac{\pi}{4}=(n+\frac{1}{2})\pi\qquad\hbox{$n$ integer}$$
The number of  zeros with imaginary part in  $[0,T]$ 
is thus $\frac{T}{
\pi}\log \frac{T}{a}-\frac{T}{ \pi}+O(1)$.

\section{Probabilistic interpretations }
We will now give interpretations of the functions 
$\xi$ and $\tilde\xi$ using first passage times of diffusions.
\subsection{Brownian motion with a drift}

The first passage time at  $x>0$ of Brownian motion started at   0
follows a 1/2 stable distribution  i.e.,
$$P(T_x\in dt)=x\frac{e^{-\frac{x^2}{ 2t}}}{ \sqrt{2\pi t^3}}dt$$
with Laplace transform
$$E[e^{-\lambda^2T_x/2}]=e^{-\lambda x}\;.$$
Adding a drift $a>0$ to the Brownian motion gives a first passage distribution
$$P^{a}(T_x\in dt)=x
\frac{e^{-\frac{x^2}{ 2t}}}
{ \sqrt{2\pi t^3}}
e^{ax-\frac{a^2t}{
2}}dt$$
with Laplace transform
$$E^{a}[e^{-\lambda^2T_x/2}]=e^{-x\sqrt{\lambda^2+a^2}+ax}.$$
This is a generalized  inverse Gaussian distribution.
In particular, its Mellin transform is
$$E^{a}[T_x^s]=(x/a)^s\frac{K_{-1/2+s}(ax)}{ K_{-1/2}(ax)}=
(x/a)^s\sqrt{\pi/ax}\,e^{ax}K_{-1/2+s}(ax)$$
which gives a probabilistic interpretation of MacDonald's function
(as a function  of $s$) as a Mellin transform of a probability distribution.
 
\subsection{Three dimensional Bessel process}
There exists a similar interpretation of the  $\xi$ function, which is
discussed in details in \cite{B}, \cite{BPY}, for example, 
Consider the first passage time at  $a>0$ of a three dimensional 
Bessel process  (i.e., the norm of a three dimensional  Brownian motion)
 starting from 0. The  Laplace transform of this hitting time
is
$$E[e^{-\frac{\lambda^2}{2}S_a}]=\frac{\lambda a}{\sinh \lambda a}\;.$$
Let  $S'_a$ be an independent copy of
$S_a$, and let
$$W_a=S_a+S_a'\;;$$
then 
the density of the distribution of  $W_a$ is obtained by inverting the Laplace
transform. One gets 
$$P(W_a\in dx)=\sum_{n=1}^\infty(\pi^4n^4x/a^4-3\pi^2 n^2/a^2)
e^{-\pi^2n^2x/2a^2}dx$$
from which one can compute the Mellin transform 
$$E[W_a^s]=2(2a^2/\pi)^{s}\xi(2s).$$
The function $2\xi$ thus has a   probabilistic interpretation, as
 Mellin transform of 
$\sqrt{\frac{\pi}{2}W_1}$. 
\subsection{Infinite divisibility}\label{infdis}
The distributions of
 $T_x$ and $W_a$ are  infinitely divisible.
Indeed
$$
\begin{array}{rcl}
\log E^a[\exp(-\frac{\lambda^2}{2}T_x)]&=& -x\sqrt{\lambda^2+a^2}+ax\\
&=&x\int_0^\infty (e^{-\frac{\lambda^2}{2}t}-1)\frac{e^{-\frac{a^2}{2}t}}{\sqrt{2\pi t^3}}dt
\end{array}
$$
which shows that $T_x$ is a subordinator with
 L\'evy measure
$$\frac{e^{-\frac{a^2}{2}t}}{\sqrt{2\pi t^3}}dt.
$$
Similarly 
$$
\begin{array}{rcl}
\log E[\exp(-\frac{\lambda^2}{2}W_a)]&=& 2\log(\lambda a/\sinh(\lambda a))\\
&=&2\int_0^\infty (e^{-\frac{\lambda^2}{2}t}-1)\sum_{n=1}^\infty
e^{-\pi^2 n^2 t/a^2}dt
\end{array}$$
therefore the variable $W_a$ has the distribution 
of a subordinator, with L\'evy measure
$$2\sum_{n=1}^\infty
e^{-\pi^2n^2 t/a^2}dt,
$$
 taken at time 1.
Observe however that the process $(W_a)_{a\geq 0}$ is not  a subordinator.
\subsection{Generalized gamma convolution }
The gamma distributions are
$$P(\gamma_{\omega,c}\in dt)=\frac{c ^{-\omega}}{\Gamma(\omega)}
t^{\omega-1}e^{- t/c}dt=\Gamma_{\omega,c}(dt)$$
where  $\omega$ and $c$ are  $>0$ parameters.
The   Laplace transform is
$$E[e^{-\lambda\gamma_{\omega,c}}]=(1+\lambda/c)^{-\omega}.$$ 
The gamma distributions form  a convolution  semigroup with respect to the
parameter
 $\omega$, i.e., 
 $$\Gamma_{\omega_1,c}*\Gamma_{\omega_2,c}=\Gamma_{\omega_1+\omega_2,c}\;.$$
The  L\'evy exponent of the gamma  semigroup is
$$\psi_{c}(\lambda)=\log(1+\lambda/c)=
\int_0^{\infty}(1-e^{-\lambda t})\frac{e^{-
ct}}{ t}dt$$
so that this is the semigroup of 
a subordinator with L\'evy measure $e^{-ct}/t\,dt$.

The generalized gamma  convolutions are the distributions of
linear combinations, with positive coefficients,
 of independent  gamma variables, and their weak limits.

One can also characterize the  generalized gamma convolutions  as
the infinitely divisible distributions with a 
L\'evy exponent of the form
$$\psi(\lambda)=\int_{0}^\infty\psi_{c}(\lambda)d\nu(c)$$
for some positive measure $\nu$ which integrates $1/c$ at $\infty$.
This measure is called the Thorin measure of the generalized gamma
 distribution.
The variables $T_x$ and $W_a$ of the preceding paragraph
are generalized gamma  convolutions.
Indeed
 it is easy to check, using the computations of section  \ref{infdis},
that $W_a$ has a generalized gamma convolution as distribution, with  Thorin measure 
\begin{equation}\label{nu11}
\nu(dc)=2\sum_{n=1}^\infty \delta_{n^2/a^2}(dc)\;;
\end{equation}
whereas
 $T_x$ is distributed as a
generalized gamma convolution with  Thorin measure
\begin{equation}\label{nu22}
\nu(dc)=\frac{dc}{
\pi\sqrt{c-a^2/2}}\qquad c>a^2/2
\end{equation}
since
$$\frac{e^{-a^2t/2}}{
\sqrt{\pi t^3}}=\int_{a^2/2}^\infty e^{- ct}\frac{dc}{
\pi\sqrt{c-a^2/2}}.$$

\subsection{Final remarks}
We can now make a connection between the preceding considerations and those of the first part of the paper. Indeed, 
 the
Thorin measures associated with the variables
$T_{x}$ and $W_{a}$ can be expressed as spectral measures associated with the generators of Brownian motion on matrix spaces. The hitting times of Brownian motion with drift are related with 
the radial part of Brownian motion in the symmetric space
$SL_2({\bf C})/SU(2)$, whereas the hitting times of the Bessel three process are related with the Brownian motion on the unitary group $SU(2)$. 
The precise relations are contained in formulas
(\ref{nu1}), (\ref{nu2}), (\ref{nu11}), (\ref{nu22}).
Thus the Riemann $\xi$ function, which is the Mellin transform of a hitting time of the Bessel three process, as in section 4.2,
and the Polya $\tilde \xi$ function from section 3.1, which appears as Mellin transform of hitting time of Brownian motion with drift,  are related in this non obvious way.

\end{document}